\documentclass[10pt,a4]{amsart}

\usepackage[mathscr]{eucal}
\usepackage{amssymb}
\usepackage{latexsym}
\usepackage{amsthm}
\usepackage{amscd}
\usepackage{amsmath}

\theoremstyle{plain}

\theoremstyle{definition}

\title{Compound basis for the space of symmetric functions}
\author{Kazuya Aokage,
Hiroshi Mizukawa 
and Hiro-Fumi Yamada}
\date{}
\address{Kazuya Aokage, Department of Mathematics, Okayama University, Okayama 700-8530, 
Japan}
\email{k.aokage@math.okayama-u.ac.jp}
\address{Hiroshi Mizukawa, Department of Mathematics,  National Defense
Academy in Japan,  Yokosuka 239-8686, Japan}
\email{mzh@nda.ac.jp}
\address{Hiro-Fumi Yamada, Department of Mathematics, Okayama University, Okayama 700-8530, 
Japan}
\email{yamada@math.okayama-u.ac.jp}
\date{}

\begin{document}
\maketitle
%\begin{flushright}{???e?????@(Kazuya Aokage)
%\\Department of Mathematics,
%\\Okayama University
%\\?????T?i?@(Hiroshi Mizukawa)
%\\Department of Mathematics,
%\\National Defense Academy
%\\?R?c?T?j?@(Hiro-Fumi Yamada)
%\\Department of Mathematics,
%\\Okayama University}
%\end{flushright}

\section{Introduction}

The aim of this note is to introduce a compound basis for the space of symmetric functions.
Our basis consists of products of Schur functions and $Q$-functions.  The basis elements are indexed by 
the partitions.  It is well known that the Schur functions form an orthonormal basis for our space.
A natural question arises.  How are these two bases connected?  In this note we present some numerical results of the transition matrix for these bases.  In particular we will see that the determinant of the transition matrix is a power of 2.  This is not a surprising fact.  However the explicit formula involves an interesting combinatorial feature.  

Our compound basis comes from the twisted homogeneous realization of the basic representation of the affine Lie algebra 
$A^{(1)}_1$ (cf. \cite{leid}).  A brief account of the relation with $A^{(1)}_1$ is found in
\cite{am}.  Also an expression of rectangular Schur functions in terms of the compound basis is given in \cite{imny}.

This note is not written in a standard style of mathematical articles.  It is more like a draft of a talk.  In particular proofs are not given here.  Details and proofs will be published elsewhere.

\vspace{10mm}

\section{Space of symmetric functions}

Throughout this note $V$ denotes the space of polynomials in infinitely many variables: 
$$V={\mathbb Q}[t_{j};j \geq 1]=\bigoplus_{n=0}^{\infty}V{(n)}.$$
Here $V(n)$ denotes the space of homogeneous polynomials of degree $n$, subject to deg $t_j =j$. This space can be regarded as the ring of symmetric functions by identifying $t_j$ with a constant multiple of the $j$-th power sum $p_j=x_1^{j}+x_2^{j}+\cdots$, where $x_k$'s are the ``original" variables. 

The first basis for $V$ consists of the Schur functions. Let $P(n)$ denote the set of the partitions of $n$. For $\lambda \in P(n)$, the Schur function $S_\lambda(t)$ indexed by $\lambda$ is defined by
 $$S_{\lambda}(t)=\sum_{\rho \in P(n)}\chi_{\rho}^{\lambda}\frac{t_{1}^{m_{1}}t_{2}^{m_{2}}\cdots}{m_{1}!m_{2}!\cdots} \quad \in V(n).$$
Here the summation runs over all $\rho=(1^{m_{1}}2^{m_{2}}\cdots) \in P(n)$, 
and the integer $\chi_{\rho}^{\lambda}$ is the irreducible character of $\lambda$ of the symmetric group $S_n$, evaluated at the conjugacy class $\rho$.
The "original" (symmetric) Schur function is recovered by putting $t_j=p_j/j.$ 
It is known that these Schur functions are orthonormal with respect to the inner product
$$\langle F,G\rangle=F(\partial)G(t)|_{t=0},$$ 
where $\partial=(\frac{\partial}{\partial t_1}, \frac{1}{2}\frac{\partial}{\partial t_2}, \frac{1}{3}\frac{\partial}{\partial t_3}, \cdots).$ 
By this orthogonality one deduces that $\{S_{\lambda}(t);\lambda \in P(n)\}$ forms an orthonormal basis for the space $V(n)$.

The second basis for $V$ is called the compound basis.
Let $SP(n)$ denote the set of the strict partitions of $n$. For $\lambda \in SP(n),$ the (Schur) $Q$-function $Q_\lambda(t)$ indexed by $\lambda$ is defined by
$$Q_{\lambda}(t)=\sum_{\rho \in OP(n)}
2^{\frac{\ell(\lambda)-\ell(\rho)+\epsilon}{2}}
\zeta_{\rho}^{\lambda}\frac{t_{1}^{m_{1}}
t_{3}^{m_{3}}\cdots}{m_{1}!m_{3}!\cdots} \quad \in V(n).$$
Here the summation runs over all odd partitions $\rho=(1^{m_{1}}3^{m_{3}}5^{m_{5}}\cdots)$ of $n$, 
the integer $\zeta_{\rho}^{\lambda}$ is the irreducible spin character of $\lambda$ of the symmetric group $S_n$, evaluated at the conjugacy class $\rho$, and  $\epsilon =0$ or $1$ according to that 
$\ell(\lambda)-\ell(\rho)$ is even or odd. 
The "original" (symmetric) $Q$-function is recovered by putting $t_j=2p_j/j.$
It is known that the $Q$-functions are orthogonal to each other with respect to the inner product
$$\langle F,G\rangle^{'}=F(2\partial)G(t)|_{t=0}$$ 
By virtue of this orthogonality, the set 
$$\{Q_{\mu}(t)S_{\nu}(t');\mu \in SP(n_0),\nu \in P(n_1),n_0+2n_1=n\}$$
 forms a basis for the space $V(n)$, where we put 
$$S_{\nu}(t')=S_{\nu}(t)|_{t_j \mapsto t_{2j}}.$$ 
This compound basis naturally arises in the study of rectangular Schur functions as weight vectors 
of the basic representation of the affine Lie algebra $A_{1}^{(1)}$ (cf. \cite{leid}, \cite{am} and \cite{imny}).

\vspace{10mm}

\section{Transition matrices}

We begin with some bijections between sets of partitions.
The first one is 
 $$\varphi:P(n) \longrightarrow  \bigcup_{n_{0}+2n_{1}=n}SP(n_{0})\times P(n_{1})$$
defined by
 $\lambda \mapsto (\lambda^{r},\lambda^{d}).$
Here the multiplicities $m_{i}(\lambda^{r})$ and
$m_{i}(\lambda^{d})$
are given respectively by
 \begin{align*}
  m_{i}(\lambda^{r})=
 \begin{cases}
 1  & m_{i}(\lambda) \equiv 1 \pmod 2\\
 0  & m_{i}(\lambda) \equiv 0 \pmod 2,
 \end{cases}
 \end{align*}
 and
  \begin{align*}
  m_{i}(\lambda^{d})=
 \begin{cases}
 \frac{1}{2}(m_{i}(\lambda)-1)  & m_{i}(\lambda) \equiv 1 \pmod 2\\
  \frac{1}{2}(m_{i}(\lambda))   & m_{i}(\lambda) \equiv 0 \pmod 2.
 \end{cases}
 \end{align*}
 For example, 
 if $\lambda=(5^34^42^71)$, then $\lambda^r=(521)$ and $\lambda^d=(54^22^3)$.
We set
 $$P({n_{0},n_{1}})=\varphi^{-1}(SP(n_{0})\times P(n_{1})).$$

The second bijection is 
$$\psi:P(n) \longrightarrow \bigcup_{n_{1}+2n_{2}=n}OP(n_{1})\times P(n_{2})$$
defined by
$\psi(\lambda)=(\lambda^o,\lambda^e)$.
Here
 $\lambda^{o}$ is obtained by 
 picking up the odd parts of $\lambda$,
 while $\lambda^e$ is obtained by 
 taking halves of the even parts. For example,  
if  $\lambda=(5^34^42^71),$ then $\lambda^o=(5^31)$ and $\lambda^e=(2^41^7)$.

The third bijection is called the 2-Glaisher map.
Let 
$\lambda=(\lambda_{1},\lambda_{2},\cdots, \lambda_{\ell})$
be a strict partition of $n$. Suppose that
$\lambda_{i}=2^{p_{i}}q_{i}\ (i=1,2,\cdots, \ell)$, where $q_{i}$ is odd.
Then an odd partition $\tilde{\lambda}$ of $n$ is defined by
$$m_{2j-1}(\tilde{\lambda})=\sum_{q_{i}=2j-1,i \geq 1}2^{p_{i}}.$$
For example, if 
$\lambda=(8,6,4,3,1)$ then $\tilde{\lambda}=(3^3,1^{13})$.
This gives a bijection between $SP(n)$ and $OP(n)$. 
%\begin{definition}
%Define a bijection $\pi$ on $P_{n}$ by 
%$$\pi(\lambda)=\psi^{-1}(\widetilde{\phi(\lambda)})\ (\lambda \in P_{n}).$$  
%Here we remark that $\pi$ gives a bijection between $\phi^{-1}(SP_{n_{1}}\times P_{n_{2}})$ and
%$\psi(OP_{n_{1}}\times P_{n_{2}})$.
%\end{definition}
%\begin{proposition}

Here are several identities for the lengths of the partitions.
Let $(n_{0},n_{1})$ be fixed. Then we have
\begin{align*}
\sum_{\lambda \in P(n)} \ell(\lambda)
 &=\sum_{\lambda \in P(n)}(\ell(\lambda^r)+2\ell(\lambda^d))
 =\sum_{\lambda \in P(n)}(\ell(\lambda^o)+\ell(\lambda^e))
 =\sum_{\lambda \in P(n)}(\ell(\tilde{\lambda}^r)+\ell(\lambda^e)),\\
 \sum_{\lambda \in P(n_{0},n_{1})}\ell(\lambda)
 &=\sum_{\lambda \in P(n_{0},n_{1})}(\ell(\lambda^r)+2\ell(\lambda^d))
 =\sum_{\lambda \in P(n_{0},n_{1})}(\ell(\lambda^o)+\ell(\lambda^e)),\\
\sum_{\lambda \in P(n)}2\ell(\lambda^d)
&=\sum_{\lambda \in P(n)}2\ell({\lambda^e})
=\sum_{\lambda \in P(n)}(\ell(\lambda^o)+\ell(\lambda^e)-\ell(\lambda^r))
=\sum_{\lambda \in P(n)}(\ell(\tilde{\lambda}^r)-\ell(\lambda^r)+\ell(\lambda^e)), \\ {\rm{and}}\\
\sum_{\lambda \in P(n_{0},n_{1})}2\ell(\lambda^d)
&=\sum_{\lambda \in P(n_{0},n_{1})}(\ell(\lambda^o)+\ell(\lambda^e)-\ell(\lambda^r)).
\end{align*}
%\end{proposition}

For simplicity we write
$$W_{\lambda}(t)=Q_{\lambda^r}(t)S_{\lambda^d}(t')$$
for $\lambda \in P(n)$.
Our problem is to determine the transition matrix between two bases. 
Let $A_{n}=(a_{\lambda \mu})$ be defined by
%\begin{equation}
$$S_{\lambda}(t)=\sum_{\mu \in P(n)}a_{\lambda \mu}W_{\mu}(t)$$
%\end{equation}
for $\lambda \in P(n)$.
Here is a small list of $A_n$.
%\begin{example}
$$A_{2}=
 \begin{array}{c|cc}
 &(2,\emptyset)&(\emptyset,1)\\
 \hline
(2)&1&1\\
 (1^2)&1&-1
 \end{array}
 $$
$$A_{3}=
 \begin{array}{c|ccc}
 &(3,\emptyset)&(21,\emptyset)&(1,1)\\
 \hline
(3)&1&0&1\\
 (21)&1&1&0\\
 (1^3)&1&0&-1
 \end{array}
 $$

 $$A_{4}=\begin{array}{c|ccccc}
&(4,\emptyset)&(31,\emptyset)&(2,1)&(\emptyset,2)&(\emptyset,1^2)\\
\hline
(4)&1&0&1&1&0\\
(31)&1&1&1&-1&0\\
(2^2)&0&1&0&1&1\\
(21^2)&1&1&-1&0&-1\\
(1^4)&1&0&-1&0&1
 \end{array}. $$
$$A_{5}=\begin{array}{c|ccccccc}
&(5,\emptyset)&(41,\emptyset)&(32,\emptyset)&(3,1)&(21,1)&(1,2)&(1,1^2)\\
\hline
(5)&1&0&0&1&0&1&0\\
(41)&1&1&0&1&1&0&0\\
(32)&0&1&1&1&0&0&1\\
(31^2)&1&1&1&0&0&-1&-1\\
(2^21)&0&1&1&-1&0&1&0\\
(21^3)&1&1&0&-1&-1&0&0\\
(1^5)&1&0&0&-1&0&0&1
 \end{array}. $$
%\end{example}

\vspace{10mm}

One readily sees that the entries are integers. Also, looking at the columns corresponding to $(\mu,\emptyset)$ with $\mu \in SP(n)$, entries are non-negative integers. The submatrix consisting of these columns will be denoted by $\Gamma_n$. The entries of $\Gamma_n $ are called the Stembridge coefficients, whose combinatorial nature has been known (\cite{st}, \cite{mac}).

Here we recall the definition of decomposition matrices for the $p$-modular representations of the symmetric group $S_n$.  Let $p$ be a fixed prime number.  A partition 
$\lambda=(\lambda_1,\lambda_2,\cdots,\lambda_{\ell})$ is said to be $p$-regular of there are no parts satisfying $\lambda_i=\lambda_{i+1}=\cdots=\lambda_{i+p-1} \geq 1.$ 
Note that a 2-regular partition is nothing but a strict partition. The set of $p$-regular partitions of $n$ is denoted by $P^{r(p)}(n).$
A partition $\rho=(1^{m_1}2^{m_2}\cdots)$ is said to be $p$-class regular of $m_{p}=m_{2p}=\cdots=0.$ Note that a 2-class regular partition is nothing but an odd partition. 
The set of $p$-class regular partitions of $n$ is denoted by $P^{c(p)}(n)$. 
The $p$-Glaisher map $\lambda \mapsto \tilde{\lambda}$ is defined in a natural way.  
This gives a bijection between $P^{r(p)}(n)$ and $P^{c(p)}(n)$. 
For $\lambda \in P^{r(p)}(n)$, we define the Brauer-Schur function $B_{\lambda}^{(p)}(t)$ indexed by 
$\lambda$ as follows.
$$B_{\lambda}^{(p)}(t)=\sum_{\rho \in P^{c(p)}(n)}\varphi_{\rho}^{\lambda}\frac{t_1^{m_1}t_2^{m_2}\cdots}{m_1!m_2!\cdots} \quad \in V(n),$$
where $\varphi_{\rho}^{\lambda}$ is the irreducible Brauer character corresponding to $\lambda$, evaluated of the $p$-regular conjugacy class $\rho$. 
These functions form a basis for the space
$V^{(p)}(n)=V^{(p)} \cap V(n)$, where 
$$V^{(p)}= {\mathbb Q}[t_j; j \geq 1, j \not \equiv 0 \, ({\rm{mod}}\, p) \}.$$

Given a Schur function $S_{\lambda}(t)$, define the $p$-reduced Schur function $S_{\lambda}^{(p)}(t)$ by "killing" all variables $t_p, t_{2p}, \cdots $;
$$S_{\lambda}^{(p)}(t)=S_{\lambda}(t)|_{t_{jp}=0}.$$ 
There $p$-reduced Schur functions are no longer 
linearly independent.  All linear relations among these polynomials are known (cf. \cite{any}).
The $p$-decomposition matrix $D_{n}^{(p)}=(d_{\lambda \mu})$ is defined by 
$$S_{\lambda}^{(p)}(t)=\sum_{\mu \in P^{r(p)}(n)}d_{\lambda \mu}B_{\mu}^{(p)}(t)$$ 
for $\lambda \in P(n)$.

Now let us go back to the case of $p=2$. We shall write $D_n$ in place of $D_{n}^{(2)}$. 
By definition, the Stembridge coefficients $a_{\lambda \mu} \quad (\lambda \in P(n),\mu \in SP(n))$ appear as 
$$S_{\lambda}^{(2)}(t)=\sum_{\mu \in SP(n)}a_{\lambda \mu}Q_{\mu}(t).$$
Looking at the matrices $D_n=(d_{\lambda \mu})$ and $\Gamma_n=(a_{\lambda \mu})$, one observes that they are 
"very similar". In fact one can prove that they are transformed to each other by column operations. 
We consider the Cartan matrix $C_n={}^tD_nD_n$ and the correspondent 
$G_n={}^t{\Gamma_n}{\Gamma_n}.$ 
The elementary divisors of $C_n$ and $G_n$ coincide. They are given by 
$\{2^{\ell(\tilde{\lambda})-\ell(\lambda)};\lambda \in SP(n)\}$ (\cite{uy}).
Our transition matrix $A_n=(a_{\lambda \mu})_{\lambda,\mu \in P(n)}$ can be regarded as a common generalization of the matrix $\Gamma_n$ of  Stembridge coefficients and 
the decomposition matrix $D_n$. 

We have a formula for the determinant of $A_n$.
$$|\det A_{n}|=2^{k_{n}},$$
where $k_{n}=\sum_{\lambda \in P(n)}\ell(\lambda^d)
=\sum_{\lambda \in P(n)}(\ell(\tilde{\lambda}^r)-\ell(\lambda^r))$.
%\begin{example}

Here is a list of $k_n$.
$$
\begin{array}{c|ccccccccc}
n&1&2&3&4&5&6&7&8&\cdots\\
\hline
k_{n}&0&1&1&4&5&11&15&28&\cdots
\end{array}$$
%\end{example}

It is natural to consider the Cartan-like matrix ${}^tA_nA_n$.  Let us look at some of these matrices.

%\begin{example}
\[
\ \ \ \ \ {}^{t}A_{2}A_2=
\bordermatrix{
  & (2,\emptyset) & (\emptyset,1)\cr
   (2,\emptyset) & 2 & 0 \cr
  (\emptyset,1) & 0 & 2 }
  \]
\[
\ \ \ \ \ {}^{t}A_{3}A_3=
\bordermatrix{
  & (3,\emptyset) & (21,\emptyset) & (1,1) \cr
   (3,\emptyset) & 3 & 1 & 0\cr
  (21,\emptyset) & 1 & 1 & 0\cr
  (1,1)& 0 & 0 &2}
  \]
  \[
  {}^{t}A_{4}A_4=
\bordermatrix{
   & (4,\emptyset)& (31,\emptyset)& (2,1) & (\emptyset, 2) & (\emptyset,1^2)\cr
  (4,\emptyset) & 4 & 2 & 0 & 0 & 0\cr
  (31,\emptyset)& 2 & 3 & 0 & 0 & 0\cr
   (2,1) & 0 & 0 & 4 & 0 & 0\cr
  (\emptyset,2) & 0 & 0 & 0 & 3 & 1\cr
  (\emptyset,1^2)& 0 & 0 & 0 & 1 & 3}
 \] 
\[
  {}^{t}A_{5}A_5=
\bordermatrix{
   &(5,\emptyset)&(41,\emptyset)&(32,\emptyset)&(3,1)&(21,1)&(1,2)&(1,1^2)\cr
  (5,\emptyset) & 5 & 3 & 1 & 0 & 0 & 0 & 0\cr
  (41,\emptyset)& 3 & 5 & 3 & 0 & 0 & 0 & 0\cr
  (32,\emptyset) & 1 & 3 & 3 & 0 & 0 & 0 & 0\cr
  (3,1)& 0 & 0 & 0 & 6 & 2 & 0 & 0\cr
  (21,1)& 0 & 0 & 0 & 2 & 2 & 0 &0\cr
  (1,2)& 0 & 0 & 0 & 0 & 0 & 3 &1\cr
  (1,1^2) & 0 & 0 & 0 & 0 & 0 & 1 & 3}
\] 
%\end{example}

\vspace{10mm}

It can be verified that ${}^tA_nA_n$ is block diagonal indexed by the pairs $(n_0,n_1)$. 
Let $B_{n_0,n_1}$ be the corresponding block in ${}^tA_nA_n$. 
Note that the principal block $B_{n,0}$ is nothing but the matrix $G_n$. 
It is plausible that there is a nice formula for elementary divisors of the block $B_{n_0,n_1}$.
At present, however, we only have a formula for the determinant.

\begin{align*}
|\det B_{n_{0},n_{1}}|
&=2^{ \sum_{\lambda \in P(n_{0},n_{1})}(\ell(\tilde{\lambda}^r)-\ell(\lambda^r)+\ell(\lambda^d))}. 
\end{align*}

\vspace{10mm}

\section{Towards the general characteristic}

We want a compound basis for the general characteristic $p$, i.e., a basis for the space
$$V= V^{(p)} \otimes V_{(p)},$$
where $V_{(p)} = \mathbb{Q}[t_{pj}; j \geq 1]$.  However, since Schur's $Q$-functions are defined only for the strict partitions, we must give up taking $Q$-functions.  We saw that the Stembridge matrix $\Gamma_n$ and the decomposition matrix $D_n$ are similar.  Therefore, for the case of general $p$, we adopt the Brauer-Schur functions $B_{\lambda}^{(p)}(t)$ instead of $Q$-functions.  

Let $p$ be a fixed prime number.  For a partition $\lambda = (\lambda_1, \cdots, \lambda_{\ell})$ of $n$, 
partitions $\lambda^{r(p)}$ and $\lambda^{d(p)}$ are defined in the following way.
The multiplicities $m_{i}(\lambda^{r(p)})$ and
$m_{i}(\lambda^{d(p)})$ are given respectively by
 \begin{align*}
  m_{i}(\lambda^{r(p)})= k \quad {\rm{if}} \quad m_{i}(\lambda) \equiv k  \pmod p
\end{align*}
 and
  \begin{align*}
  m_{i}(\lambda^{d(p)})=
 \frac{m_{i}(\lambda)-k}{p}  \quad {\rm{if}} \quad m_{i}(\lambda) \equiv k \pmod p.
 \end{align*}
For example, if $p=3$ and $\lambda = (5^34^42^{11}1^2)$, then $\lambda^{r(p)}=(42^21^2)$ and $\lambda^{d(p)}=(542^3)$.  This gives a bijection 
$$\varphi^{(p)}: P(n) \longrightarrow \bigcup_{n_0+pn_1=n}P^{r(p)}(n_0) \times P(n_1).$$
In view of this bijection, we define, for  $\lambda \in P(n)$, 
$$W_{\lambda}^{(p)}(t) = B^{(p)}_{\lambda^{r(p)}}(t) S_{\lambda^{d(p)}}(t_{(p)}),$$
where $t_{(p)}=(t_p, t_{2p},t_{3p}, \cdots).$  These functions are linearly independent and form a compound basis for the space $V(n)$.  Inconsistently $W_{\lambda}^{(2)}(t)$ and $W_{\lambda}(t)$ are not the same.  Unfortunately we do not know any connection with representation theory of affine Lie algebras yet.  At present we only verify some numerical properties of this compound basis.

Let $A_n^{(p)} = (a_{\lambda\mu})$ be the transition matrix defined by 
$$S_{\lambda}(t) = \sum_{\mu \in P(n)} a_{\lambda\mu}W_{\mu}^{(p)}(t)$$
for $\lambda \in P(n)$.  One verifies that $A_n^{(p)}$ is an integral matrix and 
$$\det A_n^{(p)} = p^{k_n^{(p)}},$$
where $k_n^{(p)}=\sum_{\lambda \in P(n)}\ell(\lambda^{d(p)})$.  As in the case of $p=2$, we consider the matrix 
$^tA_n^{(p)}A_n^{(p)}$.  This is a block diagonal matrix indexed by the pairs $(n_0,n_1)$.  Let $B_{n_0,n_1}^{(p)}$ be the corresponding block in $^tA_n^{(p)}A_n^{(p)}$.  It is obvious that the principal block $B_{n,0}^{(p)}$ coincides with the Cartan matrix $C_n^{(p)}$ at characteristic $p$.  The elementary divisors of $B_{n,0}^{(p)}=C_n^{(p)}$ are given (\cite{uy}) by
$$\{ p^{\frac{\ell(\tilde{\lambda})-\ell(\lambda)}{p-1}}; \lambda \in P^{r(p)}(n) \} ,$$
where $\tilde{\lambda}$ denotes the image of the $p$-regular partition $\lambda$ via the $p$-Glaisher map.  
For the general block we are only aware of the determinant.
$$\det B_{n_0,n_1}^{(p)} = p^{\Delta_{n_0,n_1}},$$
where 
$$\Delta_{n_0,n_1} = \sum_{(\mu, \nu) \in P^{r(p)}(n_0) \times P(n_1)} \left(\frac{\ell(\tilde{\mu})-\ell(\mu)}{p-1} + \ell(\nu)\right).$$

We give the matrices $A_n^{(p)}$ and $^tA_n^{(p)}A_n^{(p)}$ for the case $p=3$ and $n=5$.

$$A_{5}^{(3)}=\begin{array}{c|ccccccc}
&(5,\emptyset)&(2^21,\emptyset)&(41,\emptyset)&(32,\emptyset)&(31^2,\emptyset)&(2,1)&(1^2,1)\\
\hline
(5)&1&0&0&0&0&1&0\\
(41)&0&0&1&0&0&0&1\\
(32)&0&0&1&1&0&0&-1\\
(31^2)&0&0&0&0&1&0&0\\
(2^21)&1&1&0&0&0&-1&0\\
(21^3)&0&1&0&0&0&1&0\\
(1^5)&0&0&0&1&0&0&1
 \end{array}. $$

\[
  {}^{t}A_{5}^{(3)}A_5^{(3)}=
\bordermatrix{
   &(5,\emptyset)&(2^21,\emptyset)&(41,\emptyset)&(32,\emptyset)&(31^2,\emptyset)&(2,1)&(1^2,1)\cr
  (5,\emptyset) & 2 & 1 & 0 & 0 & 0 & 0 & 0\cr
  (2^21,\emptyset) & 1 & 2 & 0 & 0 & 0 & 0 & 0\cr
  (41,\emptyset)& 0 & 0 & 2 & 1 & 0 & 0 & 0\cr
  (32,\emptyset) & 0 & 0 & 1 & 2 & 0 & 0 & 0\cr
  (31^2,\emptyset)& 0 & 0 & 0 & 0 & 1 & 0 & 0\cr
  (2,1)& 0 & 0 & 0 & 0 & 0 & 3 &0\cr
  (1^2,1)& 0 & 0 & 0 & 0 & 0 & 0 & 3}
\] 

\vspace{10mm}

There may be another natural generalization using the Hall-Littlewood symmetric functions 
at root of unity (cf. \cite{mac}).
The Brauer-Schur function $B^{(p)}_{\lambda}(t)$ can be replaced by $P_{\lambda}(t;\exp{2\pi\sqrt{-1}/r})$ for any natural number $r$.   
This version of compound bases should be investigated separately.

\vspace{10mm}

\end{document}